\documentclass[reqno]{amsart}
\usepackage{graphicx}
\newtheorem{thm}{Theorem}[section]

\theoremstyle{definition}

\theoremstyle{remark}
\newtheorem{rem}[thm]{Remark}
\numberwithin{equation}{section}
\begin{document}

\title[Favard, Baxter, Gerominus, Rakhmanov, and Szeg\"o Theorems]{Favard, Baxter, Geronimus, Rakhmanov, Szeg\"o and the strong Szeg\"o theorems for
orthogonal trigonometric polynomials}%
\author{Zhihua Du}%
\address{Department of Mathematics, Jinan University, Guangzhou 510632, China}%
\email{tzhdu@jnu.edu.cn}%

\thanks{This work was initiated in Institute of
Mathematics, Free University Berlin when the author visited and
studied there from April, 2007 to April, 2008 by the support of
State Scholarship Fund of China.
The author appreciates Prof. Dr. Heinrich Begehr for his supervising and help and Prof. Dr. Hua Liu for his discussion.}%

%\date{}%
%\dedicatory{}%
%\commby{}%
% ----------------------------------------------------------------
\begin{abstract}
In this paper, we obtain some analogs of Favard, Baxter, Geronimus,
Rakhmanov, Szeg\"o and the strong Szeg\"o theorems appeared in the
theory of orthogonal polynomials on the unit circle (OPUC) for
orthogonal trigonometric polynomials (OTP). The key tool is the
mutual representation theorem for OPUC and OTP.
\end{abstract}
\maketitle
% ----------------------------------------------------------------
\section{Introduction and Preliminaries}

In \cite{dd08}, to investigate Riemann-Hilbert analysis, the author
and his collaborator established a mutual representation theorem for
OPUC and OTP which relates the isolated two classes of orthogonal
polynomials and gave four-term recurrence, Christoffel-Darboux
formula and some properties of zeros for OTP as its applications. In
fact, by the mutual representation theorem and the theory of OPUC,
we can get more analogous results appeared in the theory of OPUC for
OTP such as Favard, Baxter, Geronimus, Rakhmanov, Szeg\"o and the
strong Szeg\"o theorems, which is the theme of the present paper. To
do so, we need to introduce some notations as in \cite{dd08}.

Let $\mathbb{D}$ be the unit disc in the complex plane,
$\partial\mathbb{D}$ be the unit circle and $\mu$ be a nontrivial
probability measure on $\partial \mathbb{D}$ (i.e., with infinity
support, nonnegative and $\mu(\partial \mathbb{D})=1$). Throughout,
by decomposition, we always write
\begin{equation}
d\mu(\tau)=w(\tau)\frac{d\tau}{2\pi i\tau}+d\mu_{s}(\tau),
\end{equation}
where $\tau\in\partial \mathbb{D}$, $w(\tau)=2\pi i\tau
d\mu_{ac}/d\tau$ and $d\mu_{s}$ is the singular part of $d\mu$.

Introduce two inner products, one is complex as follows
\begin{equation}
\langle f, g\rangle_{\mathbb{C}}=\int_{\partial
\mathbb{D}}\overline{f(\tau)}g(\tau)d\mu(\tau)
\end{equation}
with norm $||f||_{\mathbb{C}}=[\int_{\partial
\mathbb{D}}|f(\tau)|^{2}d\mu(\tau)]^{1/2}$, where $f, g$ are complex
integrable functions on $\partial \mathbb{D}$. The other is
\begin{equation}
\langle f, g\rangle_{\mathbb{R}}=\int_{\partial
\mathbb{D}}f(\tau)g(\tau)d\mu(\tau),
\end{equation}
with norm $||f||_{\mathbb{R}}=[\int_{\partial
\mathbb{D}}|f(\tau)|^{2}d\mu(\tau)]^{1/2}$, where $f, g$ are real
integrable functions on $\partial \mathbb{D}$.

By the complex inner product (1.2), applying Gram-Schmidt procedure
to the following system
\begin{equation*}
\{1,z,z^{2},\ldots,z^{n},\ldots\},
\end{equation*}
where $z\in \mathbb{C}$, we get the unique system $\{\Phi_{n}(z)\}$
of monic orthogonal polynomials on the unit circle with respect to
$\mu$ satisfying
\begin{equation}
\langle\Phi_{n},\Phi_{m}\rangle_{\mathbb{C}}=\kappa_{n}^{-2}\delta_{nm}\,\,\,\text{with}\,\,\,
\kappa_{n}>0.
\end{equation}
Then the orthonormal polynomials $\varphi_{n}(z)$ on the unit circle
satisfy
\begin{equation}
\langle\varphi_{n},\varphi_{m}\rangle_{\mathbb{C}}=\delta_{nm}\,\,\,\text{and}\,\,\,
\varphi_{n}(z)=\kappa_{n}\Phi_{n}(z).
\end{equation}

For any polynomial $Q_{n}$ of $n$ order, its reversed polynomial
$Q_{n}^{*}$ is defined by
\begin{equation}
Q_{n}^{*}(z)=z^{n}\overline{Q_{n}(1/\overline{z})}.
\end{equation}

One famous property of OPUC is Szeg\"o recurrence \cite{sze}, i.e.,
\begin{equation}
\Phi_{n+1}(z)=z\Phi_{n}(z)-\overline{\alpha}_{n}\Phi^{*}_{n}(z),
\end{equation}
where $\alpha_{n}=-\overline{\Phi_{n+1}(0)}$ are called Verblunsky
coefficients. It is well known that $\alpha_{n}\in \mathbb{D}$ for
$n\in \mathbb{N}\cup\{0\}$. By convention, $\alpha_{-1}=-1$ (see
\cite{sim1}). Szeg\"o recurrence (1.7) is extremely useful in the
theory of OPUC. Especially, Verblunsky coefficients play an
important role in many problems for OPUC.

Using the real inner product (1.3) and Gram-Schmidt procedure to the
following over $\mathbb{R}$ linearly independent ordered set
  \begin{equation}
\Big\{1, \frac {z-z^{-1}}{2i}, \frac {z+z^{-1}}{2},  \ldots,
  \frac{z^{n}-z^{-n}}{2i}, \frac{z^{n}+z^{-n}}{2},  \ldots  \Big\},
  \end{equation}
where $z\in \mathbb{C}\setminus\{0\}$, we get the unique system
\begin{equation}
\{1, b_{1}\pi_{1}(z), a_{1}\sigma_{1}(z), \ldots, b_{n}\pi_{n}(z),
a_{n}\sigma_{n}(z), \ldots\}
\end{equation} of
the first ``monic" orthogonal Laurent polynomials on the unit circle
with respect to $\mu$ fulfilling (see \cite{dd08})
\begin{equation}
\langle\pi_{m},\sigma_{n}\rangle_{\mathbb{R}}=0,
\langle\pi_{m},\pi_{n}\rangle_{\mathbb{R}}=\langle\sigma_{m},\sigma_{n}\rangle_{\mathbb{R}}=\delta_{mn},\,\,\,m,n=1,2,\ldots
\end{equation}
and
\begin{equation}
a_{n}\sigma_{n}(z)=\frac{z^{n}+z^{-n}}{2}-\beta_{n}b_{n}\pi_{n}(z)-\imath_{n}a_{n-1}\sigma_{n-1}(z)
-\jmath_{n}b_{n-1}\pi_{n-1}(z)+\text{lower order}
\end{equation}
as well as
\begin{equation}
b_{n}\pi_{n}(z)=\frac{z^{n}-z^{-n}}{2i}-\varsigma_{n}a_{n-1}\sigma_{n-1}(z)
-\zeta_{n}b_{n-1}\pi_{n-1}(z)+\text{lower order},
\end{equation}
where $a_{n},b_{n}>0$, which are respectively the norms of the
``monic" orthonormal Laurent polynomials given by right hand sides
of (1.11) and (1.12),
\begin{equation}
\beta_{n}=\langle\frac{z^{n}+z^{-n}}{2},b_{n}^{-1}\pi_{n}\rangle_{\mathbb{R}},
\end{equation}
\begin{equation}
\imath_{n}=\langle\frac{z^{n}+z^{-n}}{2},a_{n-1}^{-1}\sigma_{n-1}\rangle_{\mathbb{R}},\,\,
\jmath_{n}=\langle\frac{z^{n}+z^{-n}}{2},b_{n-1}^{-1}\pi_{n-1}\rangle_{\mathbb{R}}
\end{equation}
and
\begin{equation}
\varsigma_{n}=\langle\frac{z^{n}-z^{-n}}{2i},a_{n-1}^{-1}\sigma_{n-1}\rangle_{\mathbb{R}},\,\,
\zeta_{n}=\langle\frac{z^{n}-z^{-n}}{2i},b_{n-1}^{-1}\pi_{n-1}\rangle_{\mathbb{R}}.
\end{equation}

Throughout, as a convention, take $\sigma_{0}=1$, $\pi_{0}=0$ and
$\beta_{0}=0$ as well as $a_{0}=b_{0}=1$.

In deed, identifying the unit circle with the interval $[0,2\pi)$
via the map $\theta\rightarrow e^{i\theta}$, we get the first
orthonormal trigonometric polynomials $\pi_{n}(\theta)$ and
$\sigma_{n}(\theta)$ for the linearly ordered trigonometric system
\begin{equation}
\{1, \sin\theta, \cos\theta, \ldots, \sin n\theta, \cos n\theta,
\ldots\}
\end{equation}
by the above process when $z=e^{i\theta},\,\theta\in [0,2\pi)$.

With above preliminaries, the relations of OTP and OPUC is stated by
the following theorem.
\begin{thm}[Mutual Representation \cite{dd08}]
Let $\mu$ be a nontrivial probability measure on the unit circle
$\partial \mathbb{D}=\{z: |z|=1\}$, $\{1, \pi_{n}, \sigma_{n}\}$ be
the unique system of the first orthonormal Laurent polynomials on
the unit circle with respect to $\mu$, and $\{\Phi_{n}\}$ be the
unique system of the monic orthogonal polynomials on the unit circle
with respect to $\mu$. Then for any $z\in \mathbb{C}$ and $n\in
\mathbb{N}$,
\begin{equation}
\Phi_{2n-1}(z)=z^{n-1}[a_{n}\sigma_{n}(z)+(\beta_{n}+i)b_{n}\pi_{n}(z)]
\end{equation}
and
\begin{equation}
\kappa^{2}_{2n}\Phi^{*}_{2n}(z)=\frac{1}{2}z^{n}[a^{-1}_{n}(1+\beta_{n}i)\sigma_{n}(z)
-ib^{-1}_{n}\pi_{n}(z)],
\end{equation}
where $\kappa_{2n}$ is the leading coefficient of the orthonormal
polynomial of order $2n$ on the unit circle with respect to $\mu$
and $\kappa_{2n}=\|\Phi_{2n}\|^{-1}_{\mathbb{C}}$, and $a_{n},
b_{n}, \beta_{n}$ are given by (1.11)-(1.13).
\end{thm}

Set
$\Lambda_{n}=-\frac{1}{2}[a_{n}^{-2}(1+\beta_{n}^{2})+b_{n}^{-2}]i$,
from (1.17) and (1.78), we obtain
\begin{thm}
\begin{equation}
a_{n}\sigma_{n}(z)=-\frac{1}{2}z^{-n}[\Lambda_{n}^{-1}b_{n}^{-2}iz\Phi_{2n-1}(z)-(1-\beta_{n}i)\Phi_{2n}^{*}(z)]
\end{equation}
and
\begin{equation}
b_{n}\pi_{n}(z)=-\frac{1}{2}z^{-n}[\Lambda_{n}^{-1}a_{n}^{-2}(1+\beta_{n}i)z\Phi_{2n-1}(z)-i\Phi_{2n}^{*}(z)]
\end{equation}
hold for $n\in \mathbb{N}$ and $z\in \mathbb{C}\setminus\{0\}$.
\end{thm}

As some consequences, we have (see \cite{dd08})
\begin{thm}
\begin{equation}
\kappa_{2n}^{2}=\frac{1}{4}[a_{n}^{-2}(1+\beta_{n}^{2})+b_{n}^{-2}]
\end{equation}
for $n\in \mathbb{N}\cup\{0\}$.
\end{thm}

\begin{thm}
\begin{equation}
\alpha_{2n-1}=\frac{1}{4}\kappa_{2n}^{-2}[b_{n}^{-2}-a_{n}^{-2}(1-\beta_{n}^{2})]-\frac{1}{2}\kappa_{2n}^{-2}a_{n}^{-2}
\beta_{n}i
\end{equation}
and
\begin{equation}
\alpha_{2n-2}=\frac{1}{2}(\imath_{n}-\zeta_{n}+\beta_{n-1}\varsigma_{n})
-\frac{1}{2}(\jmath_{n}+\varsigma_{n}-\beta_{n-1}\imath_{n})i
\end{equation}
for $n\in \mathbb{N}$.
\end{thm}
\begin{proof}
(1.22) is referred to \cite{dd08}. (1.23) follows from (1.11),
(1.12), (1.17) and the fact
$\alpha_{2n-2}=-\overline{\Phi_{2n-1}(0)}.$
\end{proof}

Since $\kappa_{n}^{2}/\kappa_{n+1}^{2}=1-|\alpha_{n}|^{2}$ for $n\in
\mathbb{N}\cup\{0\}$, by Theorem 1.3 and 1.4, we get

\begin{thm}
\begin{equation}
\kappa_{2n-1}^{2}=[a_{n}^{2}+b_{n}^{2}(1+\beta^{2}_{n})]^{-1}
\end{equation}
for $n\in \mathbb{N}$.
\end{thm}

Therefore, by (1.21) and (1.24), we obtain
\begin{thm}
\begin{equation}
\lim_{n\rightarrow
\infty}a_{n}b_{n}=\frac{1}{2}\exp\Big(\frac{1}{2\pi i}\int_{\partial
\mathbb{D}}\log w(\tau)\frac{d\tau}{\tau}\Big)
\end{equation}
and
\begin{equation}
\lim_{n\rightarrow
\infty}[a_{n}^{2}+b_{n}^{2}(1+\beta^{2}_{n})]=\exp\Big(\frac{1}{2\pi
i}\int_{\partial \mathbb{D}}\log w(\tau)\frac{d\tau}{\tau}\Big).
\end{equation}
\end{thm}
\begin{proof}
Note that \cite{sze,sim1}
\begin{equation}
\lim_{n\rightarrow \infty}\kappa_{n}^{-2}=\exp\Big(\frac{1}{2\pi
i}\int_{\partial \mathbb{D}}\log w(\tau)\frac{d\tau}{\tau}\Big),
\end{equation}
then (1.26) follows from (1.24) while (1.25) follows from
\begin{equation}
\kappa_{2n-1}^{2}\kappa_{2n}^{2}=\frac{1}{4}a_{n}^{-2}b_{n}^{-2}
\end{equation}
by (1.21) and (1.24).
\end{proof}

In addition, we also have
\begin{thm}
\begin{align}
&[a_{n}^{-2}(1+\beta_{n}^{2})+b_{n}^{-2}][a_{n+1}^{2}+b_{n+1}^{2}(1+\beta_{n+1}^{2})]\nonumber\\
&+(\imath_{n+1}-\zeta_{n+1}+\beta_{n}\varsigma_{n+1})^{2}+(\jmath_{n+1}+\varsigma_{n+1}-\beta_{n}\imath_{n+1})^{2}=4
\end{align}
for $n\in \mathbb{N}\cup\{0\}$.
\end{thm}
\begin{proof}
It immediately follows from (1.21), (1.23) and (1.24) since
$\kappa_{2n}^{2}/\kappa_{2n+1}^{2}=1-|\alpha_{2n}|^{2}$ for $n\in
\mathbb{N}\cup\{0\}$.
\end{proof}

\section{Analogous Theorems for Orthogonal Trigonometric Polynomials }

In present section, some analogous theorems appeared in the theory
of OPUC for orthogonal trigonometric polynomials are discussed such
as Favard, Baxter, Geronimus, Rakhmanov theorems and so on
\cite{sim1,sim2}.

\subsection{Favard Theorem} We begin with an OTP version of Favard Theorem.
Favard theorem for orthogonal polynomials on the real line is about
the orthogonality of a system of polynomials which satisfies a
three-term recurrence with appropriate coefficients \cite{enzg,man}.
Its version of OPUC is called Verblunsky theorem and well known
\cite{sim1}, that is, if $\{\alpha_{n}^{(0)}\}_{n=0}^{\infty}$ is a
sequence of complex numbers in $\mathbb{D}$, then there exists a
unique measure $d\mu$ such that $\alpha_{n}(d\mu)=\alpha_{n}^{(0)}$,
where $\alpha_{n}(d\mu)$ are the associated Verblunsky coefficients
of $d\mu$.

For orthogonal trigonometric polynomials, by Verlunsky theorem,
Theorem 1.3, 1.4 and 1.7, we have
\begin{thm}[Favard theorem for OTP] Let
$\{(a_{n}^{(0)},b_{n}^{(0)},\beta_{n}^{(0)})\}_{n=0}^{\infty}$ with
$a_{0}^{(0)},b_{0}^{(0)}=1$ and $\beta_{0}^{(0)}=0$ be a system of
three-tuples of real numbers satisfying
\begin{align}
&[(a_{n}^{(0)})^{2}+(b_{n}^{(0)})^{2}(1+(\beta_{n}^{(0)})^{2})]
[(a_{n+1}^{(0)})^{2}+(b_{n+1}^{(0)})^{2}(1+(\beta_{n+1}^{(0)})^{2})]<4(a_{n}^{(0)})^{2}(b_{n}^{(0)})^{2}
\end{align}
with $a_{n}^{(0)},b_{n}^{(0)}>0$ for $n\in \mathbb{N}\cup\{0\}$,
then there exists a nontrivial probability measure $d\mu$ on
$\partial \mathbb{D}$ such that $a_{n}(d\mu)=a_{n}^{(0)}$,
$b_{n}(d\mu)=b_{n}^{(0)}$ and $\beta_{n}(d\mu)=\beta_{n}^{(0)}$,
where $a_{n}(d\mu),b_{n}(d\mu),\beta_{n}(d\mu)$ are associated
coefficients of $d\mu$ defined by (1.11)-(1.13).
\end{thm}

\begin{proof}
For $n\in \mathbb{N}\cup\{0\}$, define \begin{equation}
\kappa_{2n}^{(0)}=\frac{1}{2}\Big[(a_{n}^{(0)})^{-2}\big(1+(\beta_{n}^{(0)})^{2}\big)+(b_{n}^{(0)})^{-2}\Big]^{\frac{1}{2}}
\end{equation}
and
\begin{equation}
\kappa_{2n+1}^{(0)}=\Big[(a_{n+1}^{(0)})^{2}+(b_{n+1}^{(0)})^{2}\big(1+(\beta_{n+1}^{(0)})^{2}\big)\Big]^{-\frac{1}{2}}.
\end{equation}

Let
\begin{equation}
\alpha_{2n-1}^{(0)}=\frac{1}{4}(\kappa_{2n}^{(0)})^{-2}\Big[(b_{n}^{(0)})^{-2}-(a_{n}^{(0)})^{-2}\big(1-(\beta_{n}^{(0)})^{2}\big)\Big]
-\frac{1}{2}(\kappa_{2n}^{(0)})^{-2}(a_{n}^{(0)})^{-2}
(\beta_{n}^{(0)})i,
\end{equation}
then $\alpha_{2n-1}^{(0)}\in\partial \mathbb{D}$ since
\begin{equation}
\Big|\alpha_{2n-1}^{(0)}\Big|^{2}=\frac{(\kappa_{2n}^{(0)})^{4}-\frac{1}{4}(a_{n}^{(0)})^{-2}(b_{n}^{(0)})^{-2}}
{(\kappa_{2n}^{(0)})^{4}}
\end{equation}
and $a_{n}^{(0)},b_{n}^{(0)}>0$. Note that (2.1) is equivalent to
\begin{equation}
\frac{\kappa_{2n}^{(0)}}{\kappa_{2n+1}^{(0)}}<1.
\end{equation}

Arbitrarily choose a sequence $\{\alpha_{2n}^{(0)}\}_{n=0}^{\infty}$
such that
\begin{equation}
\Big|\alpha_{2n}^{(0)}\Big|=\sqrt{1-(\kappa_{2n}^{(0)})^{2}\big/(\kappa_{2n+1}^{(0)})^{2}}
\end{equation}
and fix it, then $\alpha_{2n}^{(0)}\in\partial \mathbb{D}$ for $n\in
\mathbb{N}\cup\{0\}$ by (2.6).

Therefore, for the fixed sequence
$\{\alpha_{n}^{(0)}\}_{n=0}^{\infty}$, by Verblunsky theorem, there
exists a unique nontrivial probability measure $d\mu$ on $\partial
\mathbb{D}$ such that
\begin{equation}
\alpha_{n}(d\mu)=\alpha_{n}^{(0)}
\end{equation} for $n\in
\mathbb{N}\cup\{0\}$. Then for $n\in \mathbb{N}\cup\{0\}$,
\begin{equation}
\kappa_{n}(d\mu)=\kappa_{n}^{(0)}
\end{equation}
since
$\kappa_{n}(d\mu)=\prod_{j=0}^{n-1}(1-|\alpha_{j}(d\mu)|^{2})^{-\frac{1}{2}}$
(see \cite{sim1}).

Suppose that $\{\Phi_{n}(d\mu,z)\}_{n=0}^{\infty}$ is the sequence
of monic orthogonal polynomials on the unit circle with respect to
$d\mu$, set
\begin{equation}
\Sigma_{n}(z)=-\frac{1}{2}z^{-n}[(\Lambda_{n}^{(0)})^{-1}(b_{n}^{(0)})^{-2}iz\Phi_{2n-1}(d\mu,z)
-(1-\beta_{n}^{(0)}i)\Phi_{2n}^{*}(d\mu,z)]
\end{equation}
and
\begin{equation}
\Pi_{n}(z)=-\frac{1}{2}z^{-n}[(\Lambda_{n}^{(0)})^{-1}(a_{n}^{(0)})^{-2}(1+\beta_{n}^{(0)}i)z\Phi_{2n-1}(d\mu,z)
-i\Phi_{2n}^{*}(d\mu,z)]
\end{equation}
for $n\in \mathbb{N}$ and $z\in \mathbb{C}\setminus\{0\}$, where
$\Lambda_{n}^{(0)}=-\frac{1}{2}\Big[(a_{n}^{(0)})^{-2}\big(1+(\beta_{n}^{(0)})^{2}\big)+(b_{n}^{(0)})^{-2}\Big]i$.
Obviously,
\begin{equation}
\Lambda_{n}^{(0)}=-2(\kappa_{2n}^{(0)})^{2}i.
\end{equation}
By Szeg\"o recurrence and (2.8),
\begin{equation}
z\Phi_{2n-1}(d\mu,z)=\Phi_{2n}(d\mu,z)+\overline{\alpha^{(0)}_{2n-1}}\Phi^{*}_{2n-1}(d\mu,z).
\end{equation}
Hence by the orthogonality of $\Phi_{n}(d\mu, z)$ and
$\Phi_{n}^{*}(d\mu, z)$ (see \cite{sim1}), we get
\begin{equation}
\langle z^{\pm j}, \Sigma_{n}\rangle_{\mathbb{R}}=\langle z^{\pm j},
\Pi_{n}\rangle_{\mathbb{R}}=0,\,\,\,\,j=0,1,\ldots,n-1.
\end{equation}
Moreover,
\begin{equation}
\langle z^{n},
\Sigma_{n}\rangle_{\mathbb{R}}=(a_{n}^{(0)})^{2}\overline{\alpha^{(0)}_{2n-1}}+\frac{1}{2}(\kappa_{2n}^{(0)})^{-2}(1-\beta_{n}^{(0)}i),
\end{equation}
\begin{equation}
\langle z^{-n}, \Sigma_{n}\rangle_{\mathbb{R}}=(a_{n}^{(0)})^{2},
\end{equation}
\begin{equation}
\langle z^{n},
\Pi_{n}\rangle_{\mathbb{R}}=(b_{n}^{(0)})^{2}(\beta_{n}^{(0)}-i)\overline{\alpha^{(0)}_{2n-1}}+\frac{1}{2}(\kappa_{2n}^{(0)})^{-2}i,
\end{equation}
and
\begin{equation}
\langle z^{-n},
\Pi_{n}\rangle_{\mathbb{R}}=(b_{n}^{(0)})^{2}(\beta_{n}^{(0)}-i)
\end{equation}
follow from (2.9), (2.12) and the fact
$||\Phi_{n}(d\mu)||_{\mathbb{R}}^{2}=||\Phi_{n}^{*}(d\mu)||_{\mathbb{R}}^{2}=[\kappa_{n}(d\mu)]^{^{-2}}$
as well as
$(\kappa_{2n-1}^{(0)})^{2}(\kappa_{2n}^{(0)})^{2}=\frac{1}{4}(a_{n}^{(0)})^{-2}(b_{n}^{(0)})^{-2}$.
By (2.4),
\begin{equation}
\overline{\alpha^{(0)}_{2n-1}}-1=-\frac{1}{2}(\kappa_{2n}^{(0)})^{-2}(a_{n}^{(0)})^{-2}(1-\beta_{n}^{(0)}i)
\end{equation}
and
\begin{equation}
\overline{\alpha^{(0)}_{2n-1}}+1=\frac{1}{2}(\kappa_{2n}^{(0)})^{-2}\Big[(a_{n}^{(0)})^{-2}(\beta_{n}^{(0)})^{2}
+(b_{n}^{(0)})^{-2}\Big]+\frac{1}{2}(\kappa_{2n}^{(0)})^{-2}(a_{n}^{(0)})^{-2}\beta_{n}^{(0)}i.
\end{equation}
So
\begin{equation}
\langle \frac{z^{n}+z^{-n}}{2},
\Sigma_{n}\rangle_{\mathbb{R}}=(a_{n}^{(0)})^{2},\,\,\,\langle
\frac{z^{n}-z^{-n}}{2i},
\Pi_{n}\rangle_{\mathbb{R}}=(b_{n}^{(0)})^{2}
\end{equation}
and
\begin{equation}
\langle \frac{z^{n}-z^{-n}}{2i}, \Sigma_{n}\rangle_{\mathbb{R}}=0
\end{equation}
as well as
\begin{equation}
\langle \frac{z^{n}+z^{-n}}{2},
\Pi_{n}\rangle_{\mathbb{R}}=(b_{n}^{(0)})^{2}\beta_{n}^{(0)}.
\end{equation}

In addition, by straight calculations, it is easy to check that the
coefficients of $z^{n}$ and $z^{-n}$ in $\Pi_{n}(z)$ are
respectively $\frac{1}{2i}$ and $-\frac{1}{2i}$ whereas both of ones
in $\Sigma_{n}(z)-\beta_{n}^{(0)}\Pi_{n}(z)$ are $\frac{1}{2}$. Note
that (2.14) and (2.22), this fact means that $\Sigma_{n}(z)$ and
$\Pi_{n}(z)$ are just the first``monic" orthogonal Laurent
polynomials on the unit circle with respect to $d\mu$, i.e.,
\begin{equation}
\Sigma_{n}(z)=a_{n}(d\mu)\sigma_{n}(d\mu,z)\,\,\,\,\,\text{and}\,\,\,\,\,\Pi_{n}(z)=b_{n}(d\mu)\pi_{n}(d\mu,z).
\end{equation}
Since
\begin{equation}
\langle a_{n}(d\mu)\sigma_{n}(d\mu),
a_{n}(d\mu)\sigma_{n}(d\mu)\rangle_{\mathbb{R}}=a_{n}^{2}(d\mu),
\end{equation}
\begin{equation}
\langle b_{n}(d\mu)\pi_{n}(d\mu),
b_{n}(d\mu)\pi_{n}(d\mu)\rangle_{\mathbb{R}}=b_{n}^{2}(d\mu)
\end{equation}
and
\begin{equation}
\langle \frac{z^{n}+z^{-n}}{2},
b_{n}(d\mu)\pi_{n}(d\mu)\rangle_{\mathbb{R}}=b_{n}^{2}(d\mu)\beta_{n}(d\mu),
\end{equation}
therefore, by (2.21) and (2.23),
\begin{equation}
a_{n}(d\mu)=a_{n}^{(0)},\,\,\,b_{n}(d\mu)=b_{n}^{(0)},\,\,\,\beta_{n}(d\mu)=\beta_{n}^{(0)}.
\end{equation}
\end{proof}

\begin{rem}
Only for the sequence of three-tuples
$(a_{n}^{(0)},b_{n}^{(0)},\beta_{n}^{(0)})$ fulfilling (2.1), to get
(2.28), the measure $d\mu$ is not unique since the sequence can
definitely determine Verblunsky coefficients with odd subscript but
ones with even subscript from the above proof.

For $n\in \mathbb{N}$, set
\begin{equation}
\imath_{n}(d\mu)=\langle\frac{z^{n}+z^{-n}}{2},(a_{n-1}^{(0)})^{-1}\sigma_{n-1}(d\mu)\rangle_{\mathbb{R}},
\end{equation}
\begin{equation}
\jmath_{n}(d\mu)=\langle\frac{z^{n}+z^{-n}}{2},(b_{n-1}^{(0)})^{-1}\pi_{n-1}(d\mu)\rangle_{\mathbb{R}},
\end{equation}
\begin{equation}
\varsigma_{n}(d\mu)=\langle\frac{z^{n}-z^{-n}}{2i},(a_{n-1}^{(0)})^{-1}\sigma_{n-1}(d\mu)\rangle_{\mathbb{R}},
\end{equation}
and
\begin{equation}
\zeta_{n}(d\mu)=\langle\frac{z^{n}-z^{-n}}{2i},(b_{n-1}^{(0)})^{-1}\pi_{n-1}(d\mu)\rangle_{\mathbb{R}}.
\end{equation}
Then the measure $d\mu$ is unique for the sequence of seven-tuples
\begin{equation}
(a_{n}^{(0)},b_{n}^{(0)},\beta_{n}^{(0)},\imath_{n}(d\mu),\jmath_{n}(d\mu),\varsigma_{n}(d\mu),\zeta_{n}(d\mu))
\end{equation}
satisfying (2.1) by Theorem 1.4 and Verblunsky theorem. Since $d\mu$
is dependent on $(a_{n}^{(0)},b_{n}^{(0)},\beta_{n}^{(0)})$ and
$\imath_{n}(d\mu),\jmath_{n}(d\mu),\varsigma_{n}(d\mu),\zeta_{n}(d\mu)$
are dependent on $d\mu$, $a_{n}^{(0)}$ and $b_{n}^{(0)}$, then the
sequence of seven-tuples (2.33) satisfying (2.1) is dependent on the
sequence of three-tuples $(a_{n}^{(0)},b_{n}^{(0)},\beta_{n}^{(0)})$
fulfilling (2.1). Considering the uniqueness of $d\mu$ for the
sequence of (2.33) with (2.1), we call that $d\mu$ is selectively
unique for the sequence
$\{(a_{n}^{(0)},b_{n}^{(0)},\beta_{n}^{(0)})\}_{n=0}^{\infty}$
satisfying (2.1) and $a_{n}^{(0)},b_{n}^{(0)}>0$ as well as
$a_{0}^{(0)},b_{0}^{(0)}=1$ and $\beta_{0}^{(0)}=0$.
\end{rem}
\subsection{Baxter Theorem}
Let
\begin{equation}
c_{n}=\int_{\partial \mathbb{D}}\overline{\tau}^{n}d\mu(\tau),
\,\,\,n\in \mathbb{N}\cup\{0\}
\end{equation}
be moments of $\mu$, Baxter theorem for OPUC states that
$\sum_{n=0}^{\infty}|\alpha_{n}|<0$ if and only if
$\sum_{n=0}^{\infty}|c_{n}|<0$ and
$d\mu(\tau)=w(\tau)\frac{d\tau}{2\pi i\tau}$ with $w(\tau)$
continuous and $\min_{\tau\in\partial \mathbb{D}}w(\tau)>0$.

For orthogonal trigonometric polynomials, we have
\begin{thm}[Baxter theorem for OTP]
Let $\mu$ be a nontrivial probability measure on $\partial
\mathbb{D}$, $a_{n}, b_{n}, \beta_{n}$ be the associated
coefficients given in (1.11)-(1.13) and $c_{n}$ be the moments of
$\mu$ defined by (2.34), then
\begin{align}
&\sum_{n=0}^{\infty}\sqrt{1-\frac{1}{4}[a_{n}^{-2}(1+\beta_{n}^{2})+b_{n}^{-2}][a_{n+1}^{2}+b_{n+1}^{2}(1+\beta_{n+1}^{2})]}
\nonumber\\
+&\sum_{n=0}^{\infty}\sqrt{\frac{a_{n}^{4}+b_{n}^{4}(1+\beta_{n}^{2})^{2}+2a_{n}^{2}b_{n}^{2}(\beta_{n}^{2}-1)}
{a_{n}^{4}+b_{n}^{4}(1+\beta_{n}^{2})^{2}+2a_{n}^{2}b_{n}^{2}(\beta_{n}^{2}+1)}}<\infty
\end{align}
is equivalent to $\sum_{n=0}^{\infty}|c_{n}|<0$ and
$d\mu(\tau)=w(\tau)\frac{d\tau}{2\pi i\tau}$ with $w(\tau)$
continuous and $\min_{\tau\in\partial \mathbb{D}}w(\tau)>0$.
\end{thm}
\begin{proof}
It immediately follows from Theorem 1.3, 1.4 and 1.7 as well as the
Baxter theorem for OPUC.
\end{proof}

\subsection{Geronimus Theorem} To discuss Geronimus theorem, it is
necessary to introduce some basic notions of Schur algorithm (see
\cite{sim1}).

An analytic function $F$ on $\mathbb{D}$ is called a Carath\'eodory
function if and only if $F(0)=1$ and $\Re F(z)>0$ on $\mathbb{D}$.
An analytic function $f$ on $\mathbb{D}$ is called a Schur function
if and only if $\sup_{z\in \mathbb{D}}|f(z)|<1$. Let
\begin{equation}
F(z)=\int_{\partial\mathbb{D}}\frac{\tau+z}{\tau-z}d\mu(\tau)
\end{equation}
be an associated Carath\'eodory function of $\mu$, then
\begin{equation}
f(z)=\frac{1}{z}\frac{F(z)-1}{F(z)+1}
\end{equation}
is a Schur function related to $\mu$.

Starting with a Schur function $f_{0}$, Schur algorithm actually
provides an approach to continuously map one Schur function to
another by a series of transforms of the form
\begin{equation}
\begin{cases}
f_{n+1}(z)=\displaystyle\frac{1}{z}\frac{f_{n}(z)-\gamma_{n}}{1-\overline{\gamma}_{n}f_{n}(z)},\\[4mm]
\gamma_{n}=f_{n}(0).
\end{cases}
\end{equation}
$f_{n}$ are called Schur iterates and $\gamma_{n}$ are called Schur
parameters associated to $f_{0}$. Due to Schur, it is well known
that there is a one to one correspondence between the set of Schur
functions which are not finite Blaschke products and the set of
sequences of $\{\gamma_{n}\}_{n=0}^{\infty}$ in $\mathbb{D}$.
Geronimus theorem for OPUC asserts that if $\mu$ is a nontrivial
probability measure on $\partial \mathbb{D}$, the Schur parameters
$\{\gamma_{n}\}_{n=0}^{\infty}$ associated to $f_{0}$ related to
$\mu$ by (2.36) and (2.37) are identical to the Verblunsky
coefficients $\{\alpha_{n}\}_{n=0}^{\infty}$.

For orthogonal trigonometric polynomials, we have
\begin{thm}[Geronimus theorem for OTP]
Let $\mu$ be a nontrivial probability measure on $\partial
\mathbb{D}$, if $\gamma_{n}$ are Schur parameters and $a_{n}$,
$b_{n}$, $\beta_{n}$, $\imath_{n}$, $\jmath_{n}$, $\varsigma_{n}$,
$\zeta_{n}$ are coefficients associated to $\mu$ defined by
(1.11)-(1.15), then
\begin{equation}
\gamma_{2n-1}=\frac{a_{n}^{2}-b_{n}^{2}(1-\beta_{n}^{2})}{a_{n}^{2}+b_{n}^{2}(1+\beta_{n}^{2})}
-\frac{2b_{n}^{2}\beta_{n}}{a_{n}^{2}+b_{n}^{2}(1+\beta_{n}^{2})}i
\end{equation}
and
\begin{equation}
\gamma_{2n-2}=\frac{1}{2}(\imath_{n}-\zeta_{n}+\beta_{n-1}\varsigma_{n})
-\frac{1}{2}(\jmath_{n}+\varsigma_{n}-\beta_{n-1}\imath_{n})i
\end{equation}
for $n\in \mathbb{N}$.
\end{thm}

\begin{proof}
It follows from Theorem 1.3 and 1.4 as well as Geronimus theorem for
OPUC.
\end{proof}

\subsection{Rakhmanov Theorem and Szeg\"o Theorem}
Let $d\mu$ have the decomposition form (1.1),
$\{\alpha_{n}\}_{n=0}^{\infty}$ be the Verblunsky coefficients of
$\mu$, Rakhmanov theorem for OPUC states that if $w(\tau)>0$ for
a.e. $\tau\in
\partial\mathbb{D}$, then $\lim_{n\rightarrow\infty}|\alpha_{n}|=0$.
From this theorem, the OTP version of Rakhmanov theorem is the
following
\begin{thm}[Rakhmanov theorem for OTP]
Let $\mu$ be a nontrivial probability measure on $\partial
\mathbb{D}$ with the decomposition form (1.1), $a_{n}, b_{n},
\beta_{n}$ be the associated coefficients of $\mu$ given in
(1.11)-(1.13). If $w(\tau)>0$ for a.e. $\tau\in \partial\mathbb{D}$,
then
\begin{equation}
\lim_{n\rightarrow\infty}\frac{a_{n}^{4}+b_{n}^{4}(1+\beta_{n}^{2})^{2}+2a_{n}^{2}b_{n}^{2}(\beta_{n}^{2}-1)}
{a_{n}^{4}+b_{n}^{4}(1+\beta_{n}^{2})^{2}+2a_{n}^{2}b_{n}^{2}(\beta_{n}^{2}+1)}=0
\end{equation}
and
\begin{equation}
\lim_{n\rightarrow\infty}\frac{1}{4}[a_{n}^{-2}(1+\beta_{n}^{2})+b_{n}^{-2}][a_{n+1}^{2}+b_{n+1}^{2}(1+\beta_{n+1}^{2})]=1.
\end{equation}
\end{thm}

\begin{proof}
Immediate from Theorem 1.3, 1.4, 1.7 and Rakhmanov theorem for OPUC.
\end{proof}
Szeg\"o theorem for OPUC shows that
\begin{equation}
\prod_{n=0}^{\infty}(1-|\alpha_{n}|^{2})=\exp\Big(\frac{1}{2\pi
i}\int_{\partial \mathbb{D}}\log w(\tau)\frac{d\tau}{\tau}\Big).
\end{equation}

 Especially,
 \begin{equation}
 \sum_{n=0}^{\infty}|\alpha_{n}|^{2}<\infty\Longleftrightarrow
\frac{1}{2\pi i}\int_{\partial \mathbb{D}}\log
w(\tau)\frac{d\tau}{\tau}>-\infty.
 \end{equation}

Its analog for OTP is
\begin{thm}[Szeg\"o theorem for OTP]
Let $\mu$ be a nontrivial probability measure on $\partial
\mathbb{D}$ with the decomposition form (1.1), $a_{n}, b_{n},
\beta_{n}$ be the associated coefficients of $\mu$ given in
(1.11)-(1.13). Then
\begin{equation}
\prod_{n=0}^{\infty}\frac{a_{n+1}^{2}+b_{n+1}^{2}(1+\beta_{n+1}^{2})}{a_{n}^{2}+b_{n}^{2}(1+\beta_{n}^{2})}=
\exp\Big(\frac{1}{2\pi i}\int_{\partial \mathbb{D}}\log
w(\tau)\frac{d\tau}{\tau}\Big).
\end{equation}
In particular,
\begin{align}
&\sum_{n=0}^{\infty}\left\{1-\frac{1}{4}[a_{n}^{-2}(1+\beta_{n}^{2})+b_{n}^{-2}][a_{n+1}^{2}+b_{n+1}^{2}(1+\beta_{n+1}^{2})]\right\}
\nonumber\\
+&\sum_{n=0}^{\infty}\frac{a_{n}^{4}+b_{n}^{4}(1+\beta_{n}^{2})^{2}+2a_{n}^{2}b_{n}^{2}(\beta_{n}^{2}-1)}
{a_{n}^{4}+b_{n}^{4}(1+\beta_{n}^{2})^{2}+2a_{n}^{2}b_{n}^{2}(\beta_{n}^{2}+1)}<\infty
\end{align}
is equivalent to $\displaystyle\frac{1}{2\pi i}\int_{\partial
\mathbb{D}}\log w(\tau)\frac{d\tau}{\tau}>-\infty$.
\end{thm}
\begin{proof}
By Theorem 1.3, 1.4, 1.7 and Szeg\"o theorem for OPUC.
\end{proof}

\subsection{The Strong Szeg\"o Theorem}Let $d\mu$ have the decomposition form
(1.1) satisfying the Szeg\"o condition
\begin{equation}
\frac{1}{2\pi i}\int_{\partial \mathbb{D}}\log
w(\tau)\frac{d\tau}{\tau}>-\infty,
\end{equation}
it is accustomed to introduce the Szeg\"o function as follows
\begin{equation}
D(z)=\exp\Big(\frac{1}{4\pi i}\int_{\partial
\mathbb{D}}\frac{\tau+z}{\tau-z}\log w(\tau)\frac{d\tau}{\tau}\Big).
\end{equation}
It is easy to get that $D(z)$ is analytic and nonvanishing in
$\mathbb{D}$, lies in the hardy space $H^{2}(\mathbb{D})$ and
$\lim_{r\uparrow1}D(r\tau)=D(\tau)$ for a.e. $\tau\in\partial
\mathbb{D}$ as well as $|D(\tau)|^{2}=w(\tau)$. Let
\begin{equation}
D(z)=\exp\Big(\frac{1}{2}\hat{L}_{0}+\sum_{n=1}^{\infty}\hat{L}_{n}z^{n}\Big),\,\,\,z\in
\mathbb{D}.
\end{equation}
Due to Ibragimov, the sharpest form of the strong Szeg\"o theorem
for OPUC (see \cite{sim1}) says that
\begin{equation}
\sum_{n=0}^{\infty}n|\alpha_{n}|^{2}<\infty\Longleftrightarrow
d\mu_{s}=0 \,\,\,\text{and}
\,\,\,\sum_{n=0}^{\infty}n|\hat{L}_{n}|^{2}<\infty.
\end{equation}

The corresponding result for OTP can be stated as follows
\begin{thm}[The strong Szeg\"o theorem for OTP]
Let $\mu$ be a nontrivial probability measure on $\partial
\mathbb{D}$ with the decomposition form (1.1) satisfying the Szeg\"o
cindition (2.47), $a_{n}, b_{n}, \beta_{n}$ be the associated
coefficients of $\mu$ given in (1.11)-(1.13), and
$\{\hat{L}_{n}\}_{n=0}^{\infty}$ be the Taylor coefficients of the
Szeg\"o function $D(z)$ at $z=0$ which are defined by (2.48) and
(2.49). Then
\begin{align}
&\sum_{n=0}^{\infty}2n\left\{1-\frac{1}{4}[a_{n}^{-2}(1+\beta_{n}^{2})+b_{n}^{-2}][a_{n+1}^{2}+b_{n+1}^{2}(1+\beta_{n+1}^{2})]\right\}
\nonumber\\
+&\sum_{n=0}^{\infty}(2n-1)\left\{\frac{a_{n}^{4}+b_{n}^{4}(1+\beta_{n}^{2})^{2}+2a_{n}^{2}b_{n}^{2}(\beta_{n}^{2}-1)}
{a_{n}^{4}+b_{n}^{4}(1+\beta_{n}^{2})^{2}+2a_{n}^{2}b_{n}^{2}(\beta_{n}^{2}+1)}\right\}<\infty
\end{align}
is equivalent to $d\mu_{s}=0$ and
$\sum_{n=0}^{\infty}n|\hat{L}_{n}|^{2}<\infty$.
\end{thm}
\begin{proof}
From Theorem 1.3, 1.4, 1.7 and the strong Szeg\"o theorem for OPUC
of the form (2.50).
\end{proof}

In the above, by the mutual representation theorem for OTP and OPUC,
we give some analogous theorems for orthogonal trigonometric
polynomials corresponding to ones for orthogonal polynomials on the
unit circle. In fact, by the mutual representation theorem, we can
obtain more corresponding results for orthogonal trigonometric
polynomials. For example, the important and useful Bernstein-Szeg\"o
measure can be expressed in terms of orthogonal trigonometric
polynomials as follows
\begin{equation*}
d\mu_{n}=
\begin{cases}
\displaystyle\frac{a_{m}^{2}+b_{m}^{2}(1+\beta^{2}_{m})}{|a_{m}\sigma_{m}(\theta)+(\beta_{m}+i)b_{m}\pi_{m}(\theta)|^{2}}
\frac{d\theta}{2\pi},\,\,\,n=2m-1,\\[3mm]
\displaystyle\frac{a_{m}^{2}b_{m}^{2}}{a_{m}^{2}+b_{m}^{2}(1+\beta^{2}_{m})}\frac{1}{|a_{m}^{-1}(\beta_{m}-i)\sigma_{m}(\theta)
-b_{m}^{-1}\pi_{m}(\theta)|^{2}} \frac{d\theta}{2\pi},\,\,\,n=2m.
\end{cases}
\end{equation*}

% ----------------------------------------------------------------
\bibliographystyle{amsplain}

\end{document}